\documentclass[sts]{imsart}

\RequirePackage{amsthm,amsmath,amsfonts,amssymb}
\RequirePackage[numbers]{natbib}
\RequirePackage[colorlinks,citecolor=blue,urlcolor=blue]{hyperref}

\startlocaldefs
\numberwithin{equation}{section}
\theoremstyle{plain}
\newtheorem{theorem}{Theorem}
\newtheorem{assumption}{Assumption}
\newtheorem{lemma}{Lemma}
\newtheorem{proposition}{Proposition}
\DeclareMathOperator{\argmin}{argmin}

\newcommand{\R}{\mathbb{R}}
\newcommand{\df}{\hat{\mathsf{df}}}
\newcommand{\s}{\mathbf{s}}
\newcommand{\eps}{\varepsilon}
\DeclareMathOperator{\sgn}{sign}
\DeclareMathOperator{\E}{\mathbb E}
\DeclareMathOperator{\Rem}{Rem}
\DeclareMathOperator{\dv}{div}
\usepackage{cleveref}
\crefname{assumption}{Assumption}{Assumptions}
\usepackage{mathtools}

\endlocaldefs

\begin{document}

\begin{frontmatter}
\title{The Lasso error is bounded iff its active set size is bounded away from n in the proportional regime}
\runtitle{Lasso error is bounded iff its sparsity is bounded away from n}

\begin{aug}
\author[A]{\fnms{Pierre C}~\snm{Bellec}\ead[label=e1]{pierre.bellec@rutgers.edu}}
\address[A]{Rutgers University\printead[presep={\ }]{e1}.}
\end{aug}

\begin{abstract}
This note develops an analysis of the Lasso \( \hat b\) in linear models
without any sparsity or L1 assumption on the true regression vector, in
the proportional regime where dimension \( p \) and sample \( n \) are of
the same order.  Under Gaussian design and covariance matrix with spectrum
bounded away from 0 and $+\infty$, it is shown that the L2 risk is
stochastically bounded if and only if the number of selected variables is
bounded away from \( n \), in the sense that
$$
   (1-\|\hat b\|_0/n)^{-1} = O_P(1)
   \Longleftrightarrow
   \|\hat b- b^*\|_2 = O_P(1)
$$
as \( n,p\to+\infty \). The right-to-left implication rules out constant risk
for dense Lasso estimates (estimates with close to $n$ active variables),
which can be used to discard tuning parameters leading to dense estimates.

We then bring back sparsity in the picture, and revisit the precise phase
transition characterizing the sparsity patterns of the true regression vector
leading to unbounded Lasso risk---or by the above equivalence to dense Lasso
estimates.  This precise phase transition was established by
\citet{miolane2018distribution,celentano2020lasso} using fixed-point equations
in an equivalent sequence model.  An alternative proof of this phase
transition is provided here using simple arguments without relying on the
fixed-point equations or the equivalent sequence model.  A modification of the
well-known Restricted Eigenvalue argument allows to extend the analysis to any
small tuning parameter of constant order, leading to a bounded risk on one
side of the phase transition. On the other side of the phase transition, it is
established the Lasso risk can be unbounded for a given sign pattern as soon
as Basis Pursuit fails to recover that sign pattern in noiseless problems.
\end{abstract}


\end{frontmatter}

\section{Introduction}
Let $n$ be the sample size and $p$ the dimension.
Consider the Lasso estimator
\begin{equation}
\hat b = \argmin_{b\in \R^p} \tfrac12 \|Xb - y\|_2^2 + \lambda\sqrt n \|b\|_1
\label{lasso}
\end{equation}
where $X\in\R^{n\times p}$ is a design matrix with
iid $N(0,\Sigma)$ rows
for some population covariance matrix $\Sigma\in\R^{p\times p}$ with spectrum bounded away from 0 and $\infty$,
and $y$ follows a linear model
$$y = Xb^* + z$$ where the noise
$z\sim N(0,\sigma^2 I_n)$
is independent of $X$.

L1 regularization and the Lasso were introduced almost 30 years ago
\cite{tibshirani1996regression,chen1994basis} and are now recognized
as a fundamental tool in compressed sensing and high-dimensional statistics
to perform statistical estimation when $p$ is larger than $n$.
we refer the reader to the general introductions in \cite{candes2014mathematics} and in the books \cite{buhlmann2011statistics,giraud2014introduction}.
Numerous analysis have appeared throughout these three decades
to explain the success and limitations of L1 regularization
when $p$ is larger than $n$.
While it would be impossible to review all of them in this short note,
let us mention two categories of results, corresponding to
how large $p$ is compared to $n$ and whether the error
$\|\hat b - b^*\|_2$ converges to 0, in which case we will say that
the Lasso is consistent.
If $p\ggg n$ (in the sense that $p/n\to+\infty$), the Lasso is known
to be consistent under the sparsity condition
\begin{equation}
\tfrac kn\log(\tfrac pk ) \to 0 \qquad \text{ where }k=\|b^*\|_0
\label{regime-large-p}
\end{equation}
provided that the
tuning parameter $\lambda$ in \eqref{lasso} is correctly chosen.
In this regime, numerous works \cite{bickel2009simultaneous,negahban2012unified,buhlmann2011statistics,dalalyan2017prediction} explain this consistency when the tuning parameter
$\lambda$ is chosen slightly larger than $\sigma\sqrt{2 \log p}$
and if $X$ satisfies restricted eigenvalue conditions \cite{bickel2009simultaneous}. The tuning parameter can be further reduced to $\sigma\sqrt{2\log\frac pk}$ \cite{lecue2015regularization_small_ball_I,bellec2016slope}, although
anything smaller by a multiplicate constant, say
$0.99\sigma\sqrt{2\log\frac pk}$, leads to terrible results when $p\ggg n$ \cite{bellec2018nb_lsb}.
Beyond sparsity assumptions on $b^*$ as in \eqref{regime-large-p},
error bounds and consistency results in the regime $p\ggg n$ are also known
to hold if $\|b^*\|_1$ is suitably bounded from above:
this is the so-called ``slow rate bound'', see for instance
\cite[Corollary 6.1]{buhlmann2011statistics}.
Consistency and error bounds assuming upper bounds on $\|b^*\|_q, q\in(0,2)$ can also be
obtained \cite{ye2010rate}.
Random design matrices $X$ with iid $N(0,\Sigma)$ rows as in the present paper also
benefit from these analysis, since restricted eigenvalue conditions
on the population covariance $\Sigma$ transfer to the random gram
matrix $X^TX/n$ with high-probability
\cite{raskutti2010restricted}; the same can be said for subgaussian rows
using the tools in \cite{dirksen2015tail,plan_vershynin_liaw2017simple}.

Another regime, the so-called proportional regime, requires the dimension $p$ to be of the same order
as $n$, e.g.,
\begin{equation}
    \label{regime}
1 \le p/n < \gamma
\end{equation}
for some constant $\gamma$
that stays fixed as $n,p\to+\infty$ simultaneously. 
Inequality \eqref{regime} implies
$\lim \frac p n$ exists and is finite
along a subsequence of regression problems, although for the
purpose of this note we simply assume \eqref{regime}
for each regression problem
in the sequence.
Results in this proportional regime significantly depart from those
in the previous paragraph as consistency typically fails:
the error $\|\hat b - b^*\|$ is of constant order, and
a positive deterministic equivalent $\alpha_*$ such that
$(\|\hat b - b^*\|_2-\alpha_*)\to^P0$  can be precisely determined
\cite{bayati2012lasso,miolane2018distribution,celentano2020lasso}
where $\to^P$ denotes convergence in probability.
Precise phase transitions appear in this regime, including
the Donoho-Tanner phase transition for the success of Basis Pursuit
in compressed sensing \cite{donoho2005neighborliness}.
The works \cite{miolane2018distribution,celentano2020lasso} later showed that
the same phase transition precisely characterizes the sparsity levels
of $b^*$ that allow for a uniformly bounded L2 risk of the Lasso
for any constant tuning parameter $\lambda$ in \eqref{lasso}.
We will come back to this in \Cref{phase} which provides an
alternative proof of this result.

\subsection*{Contributions}
This note assumes the proportional regime throughout, with
$1 \le p/n \le \gamma$ for some constant $\gamma$ that stays fixed
as $n,p\to+\infty$, similarly to the works \cite{bayati2012lasso,miolane2018distribution,celentano2020lasso} discussed in the previous paragraph.
Throughout, we assume that the tuning parameter $\lambda$ in \eqref{lasso} and
the noise level $\sigma$ are fixed constants, and that the population covariance $\Sigma$
has bounded spectrum in the sense that its eigenvalues
belong to $[\kappa^{-1}, \kappa]$ for some constant $\kappa$.
With these normalizations, 
from the results
\cite{bayati2012lasso,miolane2018distribution,celentano2020lasso}
in this proportional regime, we expect the L2 error $\|\hat b - b^*\|$
to be of constant order.

The main contribution of this note is to highlight a surprising
behavior of the Lasso in this proportional regime, 
\emph{free of
any assumption on the true regression vector $b^*\in\R^p$}.
\Cref{theorem_main} that will be stated shortly shows that
as $n,p\to+\infty$, for any sequence of regression problems,
the Lasso risk $\|\hat b - b^*\|$ is bounded
if and only if the number of selected variables $\|\hat b\|_0$
stays bounded away from $n$, in the sense that
\begin{equation}
   \label{eq_main}
   (1-\tfrac{\|\hat b\|_0}{n})^{-1} = O_P(1)
   \Longleftrightarrow
   \|\hat b- b^*\|_2 = O_P(1).
\end{equation}
As we will discuss more deeply in the next section,
this is surprising for at least two reasons.

First, this result is free of any assumption on $b^*$:
no sparsity is assumed, and no upper bound is assumed on
$\|b^*\|_1$ or other norms: in the proportional regime
\eqref{eq_main} holds for any $b^*$ of growing dimension
as $n,p\to+\infty$.
This contrasts with most of the literature
mentioned above,
with sparsity assumptions on $b^*$, or weak sparsity assumptions
(such as bounds on $\|b^*\|_q$ for some $q\in(0,2)$).

Second,
we know that constant error bounds on $\|\hat b^R - b^*\|$
can be achieved where $\hat b^R$ is the Ridge regression estimate
or the minimum norm least-squares solution
\cite{belkin2020two,hastie2019surprises}.
Since $\hat b^R$ is dense, constant error bounds with dense estimators
are possible.
The above implication $\|\hat b - b^*\|_2 = O_P(1)
\Rightarrow (1-\|\hat b\|_0/n)^{-1} = O_P(1)$ is thus unexpected
particularity surprising.
A situation of a dense Lasso $\hat b$ (dense in the sense $\|\hat b\|_0=n$
or $1 - \|\hat b\|_0/n \lll 0$)
achieving a constant error is formally ruled out
by the above equivalence, whereas dense estimates such as Ridge regression
have no issues achieving constant error bounds with more than $n$
active variables.
This has practical application: tuning parameters $\lambda$
leading to $\|\hat b\|_0/n\approx 1$ can be discarded, as they lead to
unbounded risk.

In \Cref{phase}, we bring back sparsity in the picture,
and discuss the phase transition established in \cite{miolane2018distribution,celentano2020lasso}
characterizing the sign
patterns $\s\in\{-1,0,1\}^p$ of the true regression vector
for which \eqref{eq_main} holds. An alternative proof of this
precise phase transition is provided.

\section{The lasso error $\|\hat b - \beta^*\|^2$ 
    is bounded if and only if
$\|\hat b\|_0$ is bounded away from $n$}

Let us now state formally the main result of the paper
and its working assumption. 

\begin{assumption}
    \label{assumption}
    Let \( \gamma\ge 1 \) and \( \kappa,\sigma,\lambda>0 \) be constants
    such that
    \begin{equation}
        1 \le p/n\le \gamma,
        \qquad
        \kappa^{-1} I_p \preceq \Sigma \preceq \kappa I_p
        \label{assum:kappa}
    \end{equation}
    in the sense of psd matrices,
    $X$ has iid $N(0,\Sigma)$ rows and $\eps\sim N(0, \sigma^2 I_n)$
    is independent of $X$.
\end{assumption}

\begin{theorem}
    \label{theorem_main}
    Let \Cref{assumption} be fulfilled.
    For any $s_0\in(0,1)$, there exists $r_0'>0$
    and an event $\Omega_n(s_0)$ with
    $\mathbb P(\Omega_n(s_0)) \to 1$
    as \( n,p\to+\infty \) while \( \lambda,\sigma,\gamma,\kappa,s_0 \) remain fixed such that, in $\Omega_n(s_0)$,
    \begin{equation}
        \|\hat b \|_0 \le (1-s_0)n
        \Rightarrow
        \|\Sigma^{1/2}(\hat b - b^*)\|_2 \le r_0'.
        \label{implies_bounded_risk}
    \end{equation}
    For any $r_0\in(0,1)$, there exists $s_0'>0$
    and an event $\Omega_n(s_0)$ with
    $\mathbb P(\Omega^n(r_0)) \to 1$
    as \( n,p\to+\infty \) while \( \lambda,\sigma,\gamma,\kappa,r_0 \) remain fixed such that, in $\Omega^n(r_0)$,
    \begin{equation}
        \|\Sigma^{1/2}(\hat b - b^*)\|_2 \le r_0
        \Rightarrow
        \|\hat b \|_0 \le (1-s_0')n.
        \label{implies_bounded_sparsity}
    \end{equation}
    Consequently, for any sequence of regression problems
    with 
    \( n,p\to+\infty \) while \( \lambda,\sigma,\gamma,\kappa \) remain fixed,
    \eqref{eq_main} holds.
\end{theorem}
Above, the \( O_P(1) \) notation may hide constants depending on \( \gamma,\kappa, \lambda\): a random variable $W$ is $O_P(1)$ if 
for any $\epsilon>0$ there exists a constant $C=C(\epsilon,\gamma,\kappa,\lambda)$
such that $\mathbb P(|W|>C)\le \epsilon$.

\begin{proof}[Proof of \Cref{theorem_main}, implication \eqref{implies_bounded_risk}]
    We first prove \eqref{implies_bounded_risk} which is simpler
    and less surprising than \eqref{implies_bounded_sparsity}:
    as we detail below, in the event \( \|\hat b\|_0\le (1-s_0)n \)
    in which the Lasso has sparsity bounded away from \( n \),
    \begin{equation}
        \|\Sigma^{1/2}(\hat b - b^*)\|_2 \le c_{s_0,\gamma} \|X (\hat b - b^*)\|_2 n^{-1/2}
        \label{isometry-lasso}
    \end{equation}
    for some constant \( c_{s_0,\gamma}>0\) depending only on \( s_0,\gamma \).
    This follows because
    \begin{equation}
        \mathbb P\Bigl(
        \inf_{
            \substack{
            t\in \R, b\in \R^p: b\ne tb^*, \\
            \|b\|_0\le (1-s_0)n, 
            }
        } 
        \frac{\|X (b - t b^*)\|_2}{\|\Sigma^{1/2}(b - t b^*)\|_2}
        \ge \frac{\sqrt n}{c_{s_0,\gamma}}
        \Bigr)\to 1
        \label{isometry}
    \end{equation}
    as a consequence of bounds on the explicit density \cite{edelman1988eigenvalues} of the smallest
    eigenvalue of a Wishart$(n)$ matrix in dimension \( 1+\lfloor(1-s_0)n\rfloor \),
    and a union bound over \( \binom{n}{\lfloor (1-s_0)n \rfloor} \)
    possible subspaces.
    This argument to obtain \eqref{isometry} is detailed 
    in
    \cite[Proposition 2.10]{blanchard2011compressed} for \( \Sigma=I_p \)
    and 
    \cite[Appendix B]{bellec_zhang2019second_poincare} where \eqref{isometry}
    is proved using a slightly modified argument to handle \( \Sigma\ne I_p \)
    and a deterministic bias $b^*$.

    Now, in the event \eqref{isometry-lasso},
    obtaining a bound on the L2 risk is straightforward:
    Using that the objective function of the Lasso at \( \hat b \) 
    is smaller than at \( b^* \) we find
    \begin{align*}
        \|X(\hat b - b^*)\|_2^2
        &\le 2\eps^TX(\hat b - b^*)
        + 2 \lambda\sqrt n (\|b^*\|_1 - \|\hat b\|_1)
        \\&\le2 \eps^TX(\hat b - b^*)
        + 2 \lambda\sqrt{\gamma \kappa} n \|\Sigma^{1/2}(\hat b - b^*)\|_2
        \\&\le 
        (2\|\eps\|_2 + 2 \lambda\sqrt{\gamma\kappa} \sqrt n c_{s_0,\gamma} ) \|X(\hat b -  b^*)\|_2
    \end{align*}
    thanks to the Cauchy-Schwarz inequality for the first term
    and \( \|b^*\|_1 - \|\hat b\|_1\le \|b^*-\hat b\|_1\le \sqrt{p\kappa} \|\Sigma^{1/2}(\hat b - b^*)\| \)
    combined with \eqref{isometry-lasso} for the second term.
    Since \( \|\eps\|_2^2/\sigma^2 \) has \( \chi^2_n \) distribution,
    \begin{equation}
        \label{event-chi2}
        \mathbb P(\|\eps\|_2\le  2\sigma \sqrt n ) \ge 1-e^{-n/2},
    \end{equation}
    by \cite[Lemma 1]{laurent2000adaptive}
    so that using 
    \eqref{isometry-lasso} combined with the previous display
    grants a constant order upper bound on 
    \( \|X(\hat b - b^*)\| n^{-1/2} \), and using
    \eqref{isometry-lasso} again grants the desired constant order upper 
    bound on \( \|\Sigma^{1/2}(\hat b - b^*)\|_2 \).
\end{proof}

The proof of \eqref{implies_bounded_sparsity} is more subtle
and the result itself is more surprising:
Without any assumption on the sparsity of \( b^* \),
it is unexpected that a bounded risk for \( \hat b \) would imply
an upper bound on the number  \( \|\hat b\|_0 \) 
of active variables.
Below, let $a_+=\max(a,0)$ be the positive part of $a$ and
for a convex regularized estimator
\begin{equation}
\hat b \in \argmin_{b\in\R^p} \|y-Xb\|_2^2/2 + g(b)
\label{hbeta-g}
\end{equation}
where $g:\R^p\to\R$ is convex,
define its degrees-of-freedom in the sense of \cite{stein1981estimation}
as 
\begin{equation}
    \label{df}
    \df = \sum_{i=1}^n \frac{\partial}{\partial y_i}(x_i^T\hat b).
\end{equation}
We refer to \cite{bellec_zhang2018second_stein} and the references 
therein for the existence of the derivatives in \eqref{df}
for almost every $y$ and concentration properties of $\df$.

\begin{lemma}
    \label{my_lemma}
    Let \Cref{assumption} be fulfilled.
    Consider a convex regularized estimator
    \eqref{hbeta-g}
    for some convex function $g:\R^p\to\R$.
    Denote $R=\|\Sigma^{1/2}(\hat b - b^*)\|^2$ for brevity.
    Then for any $\mu\in(0,1]$,
    $$
    \E\Bigl[\Bigl(
            \frac{ \frac{\|X\hat b - y\|_2}{\sqrt n}}{(\sigma^2 + R)^{1/2}}
        -
        \Bigl(1-\frac{\df}{n}\Bigr)
        - \frac{\sqrt \mu R}{2\sigma^2}
        \Bigr)_+\Bigr]
        \le
    3 \sqrt \mu + 
    \frac{C \gamma^{5/4}}{ \mu n^{1/4}}
    $$
    for some absolute constant $C>0$.
\end{lemma}

The proof of \Cref{my_lemma} is given in \Cref{sec:lemma-proof}.
Below, we will take $\mu=n^{-1/8}$ so that if $R$ is bounded as in $R=O_P(1)$,
the term $\sqrt \mu R/\sigma^2$ is negligible (it converges to 0 in
probability), and the upper bound 
$3 \sqrt \mu + 
    \frac{C \gamma^{5/4}}{ \mu n^{1/4}}$
also converges to 0.  \Cref{my_lemma} thus gives us the implication
\begin{equation}
R=O_P(1)
\Rightarrow
\frac{\frac{\|X\hat b - y\|_2}{\sqrt n}}{(\sigma^2 + R)^{1/2}}
\le
\Bigl(1-\frac{\df}{n}\Bigr)
+ o_P(1).
\label{consequence_lemma}
\end{equation}
The main result of \cite{bellec2020out_of_sample} for the square loss
is that under assumptions such as strong convexity, or if the Lasso
has small enough sparsity with high-probability, the approximation
$$
\frac{\|X\hat b - y\|_2/\sqrt n}{(\sigma^2 + R)^{1/2}}
=
\Bigl(1-\frac{\df}{n}\Bigr)
+ o_P(1).
$$
is granted. This result is sometimes referred to as 
the consistency of Generalized Cross-Validation (GCV) for the generalization
error $\sigma^2+R$, where GCV is the estimator
$\|X\hat b - y\|^2/(n(1-\df/n)^2)$.
\Cref{my_lemma} is an one-sided extension of this result
under the only assumption that the risk $R=\|\Sigma^{1/2}(\hat b - b^*)\|^2$
is stochastically bounded, that is, $R=O_P(1)$.
It is one-sided because only the inequality 
\eqref{consequence_lemma} holds, and not the reverse.
It is unclear at this point if the reverse inequality holds in general
only under the assumption that $R=O_P(1)$.

To prove \eqref{implies_bounded_sparsity}, we need to bound
from below $(1-\df/n)$ in the case of the Lasso, where $\df=\|\hat b\|_0$
\cite{tibshirani2012}.
The starting point is \Cref{my_lemma} and the one-sided
inequality in \eqref{consequence_lemma} is that it allows us 
to control from below $(1-\df/n)$ with no a-priori assumption
on $\hat b$ other than its risk $R=\|\Sigma^{1/2}(\hat b - b^*)\|^2$ is
bounded.
The second ingredient is the lower bound 
$$
n \lambda^2 \|\hat b\|_0 \lesssim \|y-X\hat b\|_2^2
$$
that we will formally prove below.
Combined with \eqref{consequence_lemma}, this
gives obtain a bound of the form
$(\|\hat b\|_0/n)^{1/2} \lesssim (1-\|\hat b\|_0/n)$
which prevents $\|\hat b\|_0/n$ from being close to $1$. We now
make this precise.

\begin{proof}[Proof of \Cref{theorem_main}, implication \eqref{implies_bounded_sparsity}]
    \eqref{hbeta-g},
    Set $g(b) = \lambda \sqrt n \|b\|_1$ in \eqref{hbeta-g}
    so that $\hat b$ is the Lasso \eqref{lasso}.
    By \Cref{my_lemma} with $\mu=n^{-1/8}$
    and Markov's inequality, 
    \begin{equation}
        \label{event-lemma1}
        \frac{\|X\hat b - y\|_2}{\sqrt n(\sigma^2 + R)^{1/2}}
        -
        \Bigl(1-\frac{\df}{n}\Bigr)
        \le n^{-1/32}
        + \frac{n^{-1/16} R}{2 \sigma^2}
    \end{equation}
    has probability tending to 1 as $n,p\to+\infty$
    as $\lambda,\sigma^3,r_0,\gamma$ are held fixed.
    Consider also the event
    \begin{equation}
        \label{event-X}
    \|X\Sigma^{-1/2}\|_{op} \le \sqrt n (2 + \sqrt \gamma)
    \end{equation}
    which has probability tending to 1 as $n,p\to+\infty$,
    see for instance \cite[Theorem II.13]{DavidsonS01}
    or \cite[Corollary 7.3.3]{vershynin2018high}.
    The KKT conditions of \( \hat b \) give
    \begin{equation}
        e_j^T X^T(y-X\hat b) \in  \{-\sqrt n \lambda, ~\sqrt n \lambda\}
    \end{equation}
    for each \( j\in [p] \) such that \( \hat b_j\ne 0 \).
    Summing the squares of the previous display for all \( j\in[p] \)
    such that \( \hat b_j\ne 0 \), 
    \begin{align*}
        n \lambda^2 \df = n \lambda^2 \|\hat b\|_0
        &
        \textstyle
        \le \|y-X\hat b\|_2^2 \|X \sum_{j=1}^p e_je_j^T X^T\|_{op}
      \\&\le \|y-X\hat b\|_2^2 (2+\sqrt \gamma)^2 \kappa n
    \end{align*}
    in the event \eqref{event-X} for the last inequality.
    In the intersection of the events
    \eqref{event-X}, \eqref{event-lemma1}
    and $\{R \le r_0^2\}$,
    we have that for $n$ large enough,
    $$
    \frac{\lambda}{\sqrt \kappa(2+\sqrt \gamma)\sqrt{\sigma^2 + r_0^2}} 
    \Bigl(\frac{\df}{n}\Bigr)^{1/2}
    -
    \Bigl(1-\frac{\df}{n}\Bigr)
    \le \frac{2}{n^{1/32}}.
    $$
    If $\df/n \le 1/2$ there is nothing to prove, and
    if $\df/n\ge 1/2$ we lower bound the first term so that
    the above display provides a constant lower bound on
    $1-\df/n$.
    Taking $\Omega^n(r_0)$ as the intersection of \eqref{event-X}
    and \eqref{event-lemma1}, the proof is complete.
\end{proof}

\section{The phase transition in \NoCaseChange{\cite{miolane2018distribution,celentano2020lasso}}}
\label{phase}

In this section, we bring sparsity of $b^*$ back into the picture,
and discuss under which conditions the number of selected variables
is bounded away from $n$, in the sense 
$$
(1-\|\hat b\|_0/n)^{-1} = O_P(1)
.
$$
By the equivalence \eqref{eq_main}, under the assumptions
of \Cref{theorem_main}, this is equivalent
to boundedness of the risk, i.e.,
$\|\hat b- b^*\|_2 = O_P(1)$.

The phase transition for the boundedness of the risk was characterized
by \citet{miolane2018distribution} for $\Sigma= I_p$, and
by \citet{celentano2020lasso} for $\Sigma$ with spectrum bounded away
from 0 and $\infty$. A formulation of this phase transition
from \cite{celentano2020lasso} is as follows.

We say below that $b^*$ has sign pattern $\s\in\{-1,0,1\}^p$ if
$b_j^* \ne 0 \Rightarrow \sgn(b_j^*) = \s_j$
for all $j\in [p]$.
For any cone $T$ of $\R^p$, define the Gaussian width
$$
\mathcal G(T,\Sigma) = 
\E_{g\sim N(0,I_p)}
        \Bigl[\sup_{h\in T:h^T\Sigma h = 1} h^T\Sigma^{1/2} g\Bigr].
$$
Define $S=\{j\in[p]: |\s_j| =1\}$,
$S^c = \{j\in[p]: \s_j = 0\}$.

\begin{theorem}[\cite{miolane2018distribution} for $\Sigma=I_p$,
    \cite{celentano2020lasso} for $\Sigma\ne I_p$]
    Let \Cref{assumption} be fulfilled.
    Let $\s\in\{-1,0,1\}^p$ and define
    \begin{equation}
        \label{K}
        K=\{h\in\R^p: \sum_{j\in S^c} |h_j| \le - \s^T h\}
        .
    \end{equation}
    If $\Delta > 0$ is a constant independent of $n,p$ such that
    \begin{equation}
        \mathcal G(K,\Sigma) / \sqrt n
        \le 1 - 2\Delta,
        \label{assum}
    \end{equation}
    then for any $b^*$ with sign pattern $\s$ it holds
    $\|\Sigma^{1/2}(\hat b - b^*)\|_2 = O_P(1)$
    as $n,p\to+\infty$.
    On the other hand, if 
    \begin{equation}
        \mathcal G(K,\Sigma) / \sqrt n
        \ge 1 + \Delta,
        \label{assum_unbounded}
    \end{equation}
    then for any arbitrarily large constant $r_0>0$,
    there exists a sequence of regression problems
    with $b^*$ having sign pattern $\s$ such that
    $\|\Sigma^{1/2}(\hat b - b^*)\|_2 = O_P(1)$ does not hold
    as $n,p\to+\infty$.
\end{theorem}

Since the constant $\Delta$ is allowed to be arbitrarily small
in \eqref{assum}-\eqref{assum_unbounded}, this establishes a sharp
phase transition depending on the sign of $\mathcal G(K,\Sigma)/\sqrt n - 1$.
For isotropic designs $\Sigma=I_p$, this phase transition
coincides with the Donoho-Tanner phase transition for
the success of Basis Pursuit.

The works \cite{miolane2018distribution,celentano2020lasso}
prove this phase transition for the boundedness of the Lasso risk
using the Convex Gaussian Min-Max Theorem (CGMT) of \cite{thrampoulidis2018precise},
by characterizing the limit in probability of the Lasso risk
through a system of two nonlinear equations defined using
an equivalent sequence model, and by showing
all quantities of interest concentrate uniformly if
the problem parameters stay bounded away from the phase transition
as in \eqref{assum}.

In the next two subsections, we provide an alternative proof of this
beautiful phase transition on the boundedness of Lasso risk,
without relying on the CGMT.
\Cref{sec:RE} shows that a restricted eigenvalue condition coupled
with a simple deterministic argument is sufficient to prove that the Lasso
risk is bounded.
\Cref{sec:BP} shows that the failure of Basis Pursuit for a given
sign pattern in noiseless compressed sensing implies the unboundedness of the
Lasso risk in a related regression problem.
As we explain next, the restricted eigenvalue value condition \eqref{RE} below
is implied by \eqref{assum} for Gaussian designs,
and the failure of Basis Pursuit is implied by \eqref{assum_unbounded},
so that the results of the next two subsections provide
a simple alternative proof of this phase transition in
\cite{miolane2018distribution,celentano2020lasso}.

\subsection{From restricted eigenvalue 
to bounded Lasso risk}
\label{sec:RE}

For the cone $K$ in \eqref{K}
and some constant $\Delta_*>0$,
consider the restricted eigenvalue
condition \cite{bickel2009simultaneous}
\begin{equation}
    \label{RE}
    \forall v \in K,
    \qquad
    \Delta_* \|v\|_2
    \le \tfrac{1}{\sqrt n}  \|Xv\|_2.
\end{equation}
Gordon's Escape Through a Mesh theorem \cite{gordon1988milman}
is sufficient to prove that assumption \eqref{assum} 
implies \eqref{RE}:
Gordon's Escape Through a Mesh theorem \cite{gordon1988milman} gives
\begin{equation}
\inf_{h\in K: h^T\Sigma h=1}\|Xh\|_2
\ge
\sqrt{n-1}
- \mathcal G(L,\Sigma)
-t 
\label{gordon}
\end{equation}
with probability at least $1-e^{-t^2/2}$.
Taking $t=\Delta\sqrt n + \sqrt{n-1} - \sqrt n$
in \eqref{gordon} combined with \eqref{assum} gives
the restricted eigenvalue condition \eqref{RE}
with exponentially large probability for Gaussian design $X$.
Gordon's result has been extensively used to develop sharp
compressed sensing recovery guarantees \cite[and references therein]{chandrasekaran2012convex} or restricted isometry/restricted eigenvalue
conditions for Gaussian designs \cite{raskutti2010restricted}.
The next proposition presents a simple deterministic argument
that the restricted eigenvalue value condition \eqref{RE} is
all that is needed to prove that the Lasso risk is bounded.
Combined with the above consequence of Gordon's result, this gives
an alternative proof of the boundedness of the Lasso risk
under \eqref{assum}.

The classical restricted eigenvalue argument of \cite{bickel2009simultaneous}
or its variants 
show that the error vector $h=\hat b - b^*$ belongs to a cone
of a similar form as \eqref{K}. This argument requires the
tuning parameter to be large enough as discussed after
\eqref{regime-large-p}.
Since we are interested in the proportional regime setting where
the tuning parameter $\lambda$ in \eqref{lasso} is an arbitrarily
small constant, the classical argument fails. The following proposition
instead shows that the perturbation $v$ of the error vector
$h$
in \eqref{v} below 
belongs to $K$ (in particular, it is not claimed throughout the 
following proof that $h\in K$).

\begin{proposition}
    \label{prop:bounded_risk}
    Let $\lambda,\Delta_*>0$ be constants.
    Let $\s\in\{-1,0,1\}$ be a sign pattern with $k$ nonzero entries.
    Assume that there exist $c_1,c_2,c_3$ such that
    \begin{equation}
        \label{c_1c_2c_3}
        \frac{\|X\s\|_2}{\sqrt{n k}} \le c_1,
    \quad
    \frac{\sqrt n}{\sqrt k} \le c_2,
    \quad
    \frac{\|\eps\|_2}{\sqrt n}\le c_3.
    \end{equation}
    Assume the restricted eigenvalue condition
    \eqref{RE}
    for the cone $K$ in \eqref{K}.
    Then for the Lasso error 
    $h=\hat b - b^*$,
    \begin{align}
        \label{pre-risk}
        \tfrac{1}{\sqrt n}\|X h\|_2
    &\le 
    \Delta_*^{-1}[
        \lambda \sqrt{k/n}
        +c_1 c_2  c_3
    ],
    \\
        \label{L2-risk}
    \|h\|_2
    &\le
    \tfrac{1}{\sqrt n}
    \|Xh\|_2
    [
    \tfrac{c_2 c_3}{\lambda}
    +
    \Delta_*^{-1}
        (
            1+
            \tfrac{c_1c_2c_3}{\lambda}
        )
    ].
    \end{align}
\end{proposition}
The constant $c_2$ can be obtained as long as the sparsity $k$
is of the same order as $n$.
If $X$ has iid centered rows with finite covariance
with \eqref{assum:kappa},
and the noise $\eps$ centered iid components with variance $\sigma^2$,
for any $\epsilon>0$
the existence of $c_1,c_3$ such that
\eqref{c_1c_2c_3} holds with probability $1-\epsilon$
follows from Markov's inequality.
Thus the assumption \eqref{c_1c_2c_3} is mild and can be satisfied
in typical settings of the proportional regime where
sparsity, sample size and dimension are of the same order.
If $\Delta_*,\lambda,c_1,c_2,c_3$ are of constant order,
\eqref{pre-risk}-\eqref{L2-risk} show that the risk of the Lasso
is bounded as desired.

\begin{proof}[Proof of \Cref{prop:bounded_risk}]
    The KKT conditions of the Lasso problem
    (see, e.g., Lemma~1 in \cite{bellec2016prediction})
    provide
    \begin{align*}
    \|X(\hat b - b^*)\|_2^2
    &\le \eps^TX(\hat b - b^*) 
    + \lambda \sqrt n\bigl( \|b^*\|_1 - \|\hat b\|_1\bigr)
  \\&\le
  \eps^TXh
  +
  \lambda \sqrt n\bigl(-\s^Th - \|h_{S^c}\|_1 \bigr) 
    \end{align*}
    where $h=\hat b - b^*$.
    The first term $\eps^TXh$ in the second line is bounded from above by
    $\|\eps\|_2 \|Xh\|_2$. We force the introduction of an inner product
    with the sign vector $\s$ using
    $$
    \|\eps\|_2 \|Xh\|_2
    =
    \lambda \sqrt n 
    \s^T
    \Bigl(
        \s
        \frac{
    \|\eps\|_2 \|Xh\|_2}{\lambda \sqrt n k}
    \Bigr)
    $$
    where $k=\|\s\|_2^2$ is the number of nonzero signs. Then
    $$
    \|Xh\|_2^2
    \le \lambda \sqrt n
    \Bigl[
    \s^T\Bigl(
    \s
    \frac{\|Xh\|_2\|\eps\|_2}{k\lambda \sqrt n} 
    - h_S\Bigr)
    - \|h_{S^c}\|_1
    \Bigr]
    $$
    Since the left-hand side is non-negative,
    the vector
    \begin{equation}
        v=
        \s \frac{\|Xh\|_2\|\eps\|_2}{k\lambda \sqrt n} - h
        \label{v}
    \end{equation}
    satisfies
    $v\in K$, hence by \eqref{RE} it holds
    $\Delta_* \sqrt n \|v\|_2 \le \|Xv\|_2$.
    Using the Cauchy-Schwarz inequality
    $\s^T v \le \sqrt k \|v\|_2$
    we find
    \begin{align*}
    \|Xh\|_2^2
    \le \frac{\lambda \sqrt k \|Xv\|_2}{\Delta_*}
    \le
    \frac{\lambda \sqrt k
    \|Xh\|_2}{\Delta_*}
    \Bigl[1
        + \frac{\|X \s\|_2 \|\eps\|_2}{\lambda k \sqrt n}
    \Bigr]
    \end{align*}
    by the triangle inequality for the second inequality.
    Simplifying both sides by $\|Xh\|_2$ gives
    \eqref{pre-risk}.

    We now prove \eqref{L2-risk} by using
    the triangle inequality, in order to transfer to
    one-sided isometry inequality \eqref{RE}
    from $v$ in \eqref{v}
    to the error vector $h=\hat b - b^*$:
    \begin{align*}
        \|h\|_2
        &\le \|\s\|_2 \|Xh\|_2 \|\eps\|_2/(k\lambda \sqrt n)
        + \|v\|_2
      \\&\le \|Xh\|_2 \|\eps\|_2 / (\sqrt k \lambda \sqrt{n })
      + \|X v\|_2/(\sqrt n\Delta_*)
      \\&\le
      \frac{\|Xh\|_2}{\sqrt n}
      \Bigl[
      \frac{\|\eps\|_2}{\sqrt k\lambda}
      + \frac{1}{\Delta_*}
      \Bigl(1
      + \frac{\|X\s\|_2\|\eps\|_2}{k\lambda\sqrt n}
      \Bigr)
      \Bigr]
    \end{align*}
    by the triangle inequality and \eqref{v}.
    This proves \eqref{L2-risk}.
\end{proof}

\subsection{From failure of basis pursuit
to unbounded risk}
\label{sec:BP}

Consider a true signal $b_0\in \R^p$ with sign pattern $\s$,
and the Basis Pursuit problem
\begin{equation}
\tilde b_0 = \argmin_{b\in\R^p} \|b\|_1 \text{ subject to } Xb=Xb_0.
\label{BP}
\end{equation}
\citet{amelunxen2014living} establish a sharp phase transition
for the success of Basis Pursuit in recover $b_0$:
if \eqref{assum} holds then $b_0$ is the unique solution
to \eqref{BP} (in other words, \eqref{BP} succeeds at recover $b_0$),
while if
\eqref{assum_unbounded} holds then there exists $\tilde b_0$ different
than $b_0$ with $\tilde b_0$ solution to \eqref{BP}
(i.e., \eqref{BP} fails to recover $b_0$).
This follows from Theorem~II of \cite{amelunxen2014living} 
for Gaussian design with covariance $\Sigma$:
The failure of \eqref{BP} to recover $b_0$ under \eqref{assum_unbounded}
is established by
showing that the kernel of $X\Sigma^{-1/2}$ has
intersection with $\Sigma^{1/2}K$ larger than $\{0\}$,
where $K$ is the closed convex 
cone in \eqref{K}. This implies that the solution set
of \eqref{BP} is not equal to the singleton $\{b_0\}$
with high-probability:
$$
\exists \tilde b_0\ne b_0 \text{ such that }\|\tilde b_0\|_1 \le \|b_0\|_1.
$$
For the application below, it is useful to obtain a strict inequality
$\|\tilde b_0\|_1 < \|b_0\|_1$ using the following argument.
The kernel of $X\Sigma^{-1/2}$ has dimension $p-n$ and rotationally invariant distribution.
Thus, with probability one, its intersection with any deterministic
hyperplane is $\{0\}$.
Since $\Sigma^{1/2}K$ is an intersection of a finite number of half-spaces,
it follows that with probability one, if $X\Sigma^{-1/2}$ has a non-trivial
intersection with $\Sigma^{1/2}K$ then it has a non-trivial intersection
with the interior of $\Sigma^{1/2}K$ and in this case
there exists $h_0\in\text{interior}(K)$ such that
$X h_0=0$, and $h_0 \in\text{interior}(K)$.
This implies $\|b_0 + t h_0\|_1 < \|b_0\|_1$ for $t>0$ small enough.
This argument shows that for Gaussian design, under \eqref{assum_unbounded},
the solution \eqref{BP} has L1 norm strictly smaller than $\|b_0\|_1$
and
\begin{equation}
\mathbb P(
\exists \tilde b_0\ne b_0: \|\tilde b_0\|_1 < \|b_0\|_1,
~ Xb_0 = X\tilde b_0
) \ge 0.9.
\label{b0}
\end{equation}
holds. The result of \cite{amelunxen2014living} is stronger and shows
that this probability converges to 1 exponentially fast; for
simplicity of exposition we use the explicit absolute constant $0.9$.

\begin{proposition}
    \label{prop_if_BP_fails}
    Assume that the sign pattern $\s$ is such that
    we can find a deterministic $b_0\in\R^p$ with sign pattern $\s$ such that
    \eqref{b0} holds.
    Then for any arbitrarily large $r_0>0$,
    there exists
    $b^*\in\R^p$ with sign pattern $\s$ such that the Lasso
    error is at least $r_0$ with positive probability in the sense
    $\mathbb P(
    \|\Sigma^{1/2}(\hat b - b^*)\|_2
    \ge r_0
    ) \ge 0.6$.
\end{proposition}
In other words, as soon as basis pursuit does not recover $b_0$
as in \eqref{b0}, we can construct a related regression problem
in which the Lasso has unbounded risk. This argument
relies on \eqref{b0} only and not in the Gaussian design assumption.
The proof argument given below is a variant of the lower bound in
\cite[Theorem 3.1]{bellec2018nb_lsb} that the Lasso risk is bounded from below
by the compatibility constant.

Together with the fact that \eqref{b0} holds under \eqref{assum_unbounded}
by the discussion preceding \eqref{b0},
this provides an alternative proof of the unbounded risk
of the Lasso for some sparse regression vector when
\eqref{assum_unbounded} holds, which was initially established in
Proposition~14 in \cite{celentano2020lasso}
as a consequence of the CGMT.

\begin{proof}[Proof of \Cref{prop_if_BP_fails}]
    Given
    this deterministic $b_0$, set as the true regression vector
    $b^*$, and define a random $\tilde b$, by
    \begin{equation}
        \label{t}
    b_* = t b_0,
    \qquad
    \tilde b = t \tilde b_0
    \end{equation}
    for some $t>0$ that will be chosen later.

    Let $\mu\in (0,1]$.
    For a convex $g:\R^p\to \R$,
    consider 
    \begin{equation}
        \label{b_mu}
        b^\mu =
        \argmin_{b\in \R^p} L_\mu(b),
    \end{equation}
    where the objective function has an extra quadratic term:
    \begin{equation*}
        L_\mu(b) \coloneqq \tfrac 1 2  \|Xb - y\|_2^2 + g(b)
        + \tfrac{\mu n }{2}\|\Sigma^{1/2}(\hat b - b^*)\|_2^2.
        \nonumber
    \end{equation*}
    Since $b^\mu$ minimizes \( L_\mu(\cdot) \)
    we have if $g(b)= \lambda\sqrt n\|b\|_1$
    \begin{align*}
        L_\mu(b^\mu)
        &=
        L_0(b^\mu) +
        (\mu n/2) \|\Sigma^{1/2}(b^\mu - b^*)\|_2^2
        \\
        &\le L_\mu(\hat b) = L_0(\hat b) + (\mu n/2) \|\Sigma^{1/2}(\hat b - b^*)\|_2^2,
    \end{align*}
    and since \( L_0(\hat b)\le L_0(b^\mu) \) we obtain
    \begin{equation}
        \label{monotonicity}
      \|\Sigma^{1/2}(b^\mu - b^*)\|_2^2
      \le
      \|\Sigma^{1/2}(\hat b - b^*)\|_2^2
        . 
    \end{equation}
    That is, the risk of \( b^\mu \) is always smaller or equal than that
    of \( \hat b \).
    Since we are trying to prove a lower bound
    on $\|\Sigma^{1/2}(\hat b - b^*)\|_2$,
    due to \eqref{monotonicity} it is sufficient to 
    obtain a lower bound on $\|\Sigma^{1/2}(b^\mu - b^*)\|_2$.
    For brevity, set $h^\mu = b^\mu-b^*$.
    The KKT conditions of $b^\mu$ read
    $$0 \in \lambda\sqrt n \partial \|b^\mu\|_1
    +
    X^T(Xh^\mu - \eps) 
    + \mu n \Sigma h^\mu
    .$$
    Multiplying these KT conditions by $\tilde b -b^\mu$,
    using that $X\tilde b = Xb^*$, we find
    \begin{multline}
    0 \le
    \eps^TX h^\mu
    - \|X h^\mu\|_2^2
    + \mu n (\tilde b - b^*) \Sigma h^\mu
    \\
    +
    \lambda \sqrt n(\|\tilde b \|_1- \|b^\mu\|_1)
    - \mu n \|\Sigma^{1/2} h^\mu\|_2^2.
    \label{to_sum}
    \end{multline}
    We have $\eps^TX h^\mu - \|X h^\mu\|_2^2
    \le
    \|\eps\|_2^2/4$
    for the first two terms in the first line
    thanks to $ab - a^2 \le b^2/4$.
    Next, by the Cauchy-Schwarz inequality
    and \eqref{assum:kappa}
    \begin{align*}
    \lambda\sqrt n (\|b^*\|_1 - \|b^\mu\|_1)
    &\le \lambda\sqrt n \sqrt p \|h^\mu\|_2
  \\&\le \lambda\sqrt n \sqrt p \|\Sigma^{1/2} h^\mu\|_2 \sqrt \kappa
  \\&\le
  \lambda^2 p \kappa/(4\mu) + \mu n \|\Sigma^{1/2} h^\mu\|_2^2
    \end{align*}
    again thanks to $ab \le b^2/4 + a^2$.
    Summing these inequalities with \eqref{to_sum} gives
    $$
    \lambda\sqrt n(\|b^* \|_1- \|\tilde b\|_1)
    \le 
    \frac{\lambda^2 p \kappa}{4\mu}
    +\frac{\|\eps\|_2^2}{4}
    + \mu n (\tilde b - b^*) \Sigma h^\mu.
    $$
    We bound the rightmost term
    by the Cauchy-Schwarz inequality 
    $(\tilde b - b^*) \Sigma h^\mu
    \le \|\Sigma^{1/2}(\tilde b - b^*)\|_2
    \|\Sigma^{1/2} h^\mu\|_2$,
    and divide by $\mu n \|\Sigma^{1/2}(\tilde b - b^*)\|_2$, which gives
    $$
    \frac{\lambda\sqrt n(\|b^* \|_1- \|\tilde b\|_1)}{\mu n \|\Sigma^{1/2}(\tilde b - b^*)\|_2}
    \le 
    \frac{\lambda^2 p \kappa + \mu\|\eps\|_2^2}{4 \mu^2n \|\Sigma^{1/2}(\tilde b - b^*)\|_2}
    + \|\Sigma^{1/2} h^\mu\|_2.
    $$
    We now find an event of positive probability such that,
    for well chosen $\mu$ and $t$ in \eqref{t},
    the left-hand side is at least $2r_0$, and the first term on the
    right-hand side is smaller than $r_0$.
    Since all vectors involved in the left-hand side
    are proportional to $t$ by \eqref{t}, the left-hand side equals
    $$
    \frac{1}{\mu} W,
    \qquad W = 
    \frac{\lambda\sqrt n(\|b_0 \|_1- \|\tilde b_0\|_1)}{n \|\Sigma^{1/2}(\tilde b_0 - b_0)\|_2}.
    $$
    By \eqref{b0},
    the random variable $W$ satisfies $\mathbb P(W > 0) \ge 0.9$.
    By continuity of the probability, there exists a positive deterministic
    $\rho>0$ such that $\mathbb P(W > \rho) \ge 0.8$.
    Set $\mu = \rho /(2r_0)$ so that
    $\mathbb P( 2r_0 \le W/\mu ) \ge 0.8$.
    Now $\mu$ is fixed, and we can still choose $t>0$.
    The first term in the right-hand side is
    $$
    \frac{1}{t} W'
    \quad
    \text{ where }
    \quad
    W' = \frac{\lambda^2 p \kappa + \mu\|\eps\|_2^2}{4\mu^2n \|\Sigma^{1/2}(\tilde b_0 - b_0)\|_2}.
    $$
    The random variable $W'$ is positive and finite with probability
    at least $0.9$ by \eqref{b0}, so we may find
    $\rho'>0$ such that $\mathbb P(W' \le \rho') \ge 0.8$.
    Choosing $t = \rho'/r_0$ and using the union bound, have
    have with probability at least $0.6$ that
    $2r_0 \le r_0 + \|\Sigma^{1/2}h^\mu\|_2$ as desired.
\end{proof}

\appendix
\section{Proof of Lemma~\NoCaseChange{\ref{my_lemma}}}
\label{sec:lemma-proof}

\begin{proof}[Proof of \Cref{my_lemma}]
    Consider again the oracle $b^\mu$ in \eqref{b_mu}
    and recall that the monotonicity property
    \eqref{monotonicity} holds.
    We now prove that \( \|y-X b^\mu\|/\sqrt n \) and
    \(  \|y-X\hat b\|/\sqrt n \) are close
    when $\mu$ is small. Since
    \( b^\mu \) minimizes both \( L_\mu \) and
    \(b\mapsto  L_\mu(b) - \|X(b^\mu-b)\|^2/2 \), we have
    \begin{align*}
        L_\mu(b^\mu)
        &\le L_\mu(\hat b) - \|X(\hat b - b^\mu)\|_2^2/2
      \\&\le L_0(\hat b) + \frac{\mu n}{2} \|\Sigma^{1/2}(\hat b -b^*)\|_2^2
        - \|X(\hat b - b^\mu)\|_2^2/2.
    \end{align*}
    We further bound the right-hand side using
    \( L_0(\hat b)\le L_0(b^\mu) \le L_\mu(b^\mu) \) and conclude that
    \begin{equation}
        \label{bound-residual}
        \|X(\hat b - b^\mu)\|_2^2 \le \mu n  \|\Sigma^{1/2}(\hat b - b^*)\|_2^2.
    \end{equation}
    By Theorem 3.1 in \cite{bellec2020out_of_sample}
    and thanks to the strong convexity term proportional to \( \mu \)
    in the objective function of \( b^\mu \),
    $$
        \Rem_0 \coloneqq
        \Bigl(1-\frac{\df^\mu}{n}\Bigr)^2
        -
        \frac{
            \|y-X b^\mu\|_2^2/n
        }{(\sigma^2 + \|\Sigma^{1/2}(b^\mu - b^*)\|_2^2)},
    $$
    where $\df^\mu$ is the degrees of freedom of \( b^\mu \)
    defined as in \eqref{df} with $\hat b$ replaced by $b^\mu$,
    satisfies
    \begin{equation}
    \E[
        |
        \Rem_0
        |
        ]
        \le C \gamma^{5/2} \mu^{-2} n^{-1/2}
        .
        \label{10}
    \end{equation}
    While the dependence on $(\gamma,\mu)$ in the right-hand side
    is not explicit in \cite{bellec2020out_of_sample},
    it is made explicit
    in \cite[Theorem 1 with $c=1$ and $m=l$]{bellec2023corrected}.
    Using $(\sqrt a - \sqrt b)^2\le |a-b|$ for positive $a,b$,
    we obtain that
    $$
    \Rem =
    \frac{\|y-X b^\mu\|_2/\sqrt n}{(\sigma^2 + \| \Sigma^{1/2}(b^\mu - b^*)\|_2^2)^{1/2}}
    -
    \Bigl(1-\frac{\df^\mu}{n}\Bigr)
    $$
    satisfies $\E[\Rem^2] \le C \gamma^{5/2} \mu^{-2} n^{-1/2}$, which implies
    \begin{equation}
    \E|\Rem| \le 
    \sqrt{C} \gamma^{5/4} \mu^{-1}
    n^{-1/4}.
    \label{bound-Rem}
    \end{equation}
    Using the monotonicity 
    \eqref{monotonicity} and \eqref{bound-residual}, we have
    \begin{align}
    \Bigl(1- \frac{\df^\mu}{n}\Bigr) + \Rem
    &= \frac{\|y-X b^\mu\|/\sqrt n}{(\sigma^2 + \| \Sigma^{1/2}(b^\mu - b^*)\|_2^2)^{1/2}}
    \nonumber
  \\&\ge 
    \frac{\|y-X b^\mu\|_2/\sqrt n}{(\sigma^2 + R)^{1/2}}
    \quad \text{ by \eqref{monotonicity}}
    \nonumber
  \\&\ge 
    \frac{\|y-X \hat b\|_2/\sqrt n}{(\sigma^2 + R)^{1/2}}
    - \frac{\sqrt \mu R^{1/2}}{(\sigma^2 + R)^{1/2}}.
    \label{combined}
    \end{align}
    It remains to explain that $\df$ and
    $\df^\mu$ are close.
    Conditionally on \( X \), the vector field
    \( f(\eps) = X(\hat b - b^\mu) \) is 2-Lipschitz
    \cite{bellec2016bounds}. Applying Proposition 6.4
    in \cite{bellec2020out_of_sample}, there exists \( Z\sim N(0,1) \)
    and two random variables \( T,\tilde T \) with \( \E[T^2]\le1,\E[\tilde T^2]\le 1 \) such that almost surely,
    \begin{multline*}
        |\sigma^2 \dv f(\eps) - \eps^Tf(\eps)|
        \\\le \sigma |Z| \|f(\eps)\| + \sigma^2 (2 |T| + |Z\tilde T| ) 2 \sqrt n.
    \end{multline*}
    Since $\dv f(\eps) = \df - \df^\mu$
    by definition of $\df$ and $\df^\mu$,
    by the triangle inequality
    and $|\eps^Tf(\eps)| \le \|\eps\|_2\|f(\eps)\|_2$,
    the difference
    $|\df - \df^\mu|$ is bounded from above by
    \begin{align*}
    &(|Z| + \frac{\|\eps\|_2}{\sigma}) \frac{\|f(\eps)\|_2}{\sigma} + (2 |T| + |Z\tilde T| ) 2 \sqrt n.
    \\&\le
    \frac{\sqrt \mu
    (|Z| + \frac{\|\eps\|_2}{\sigma})^2}{2}
    + \frac{\|f(\eps)\|_2^2}{2\sqrt \mu\sigma^2}
    + (2 |T| + |Z\tilde T| ) 2 \sqrt n 
    \end{align*}
    thanks to
    $ab\le \frac{\sqrt \mu a^2}{2} + \frac{b^2}{\sqrt \mu 2}$.
    We use again \eqref{bound-residual} to bound
    $\|f(\eps)\|_2^2\le \mu n R$ so that, still
    almost surely,
    \begin{equation}
        \label{combined2}
        \frac{|\df - \df^\mu|}{n}
        - 
        \frac{\sqrt \mu R}{2\sigma^2}
    \le
    \Rem'
    \end{equation}
    where $\Rem'=
    \frac{\sqrt \mu
    (|Z| + \frac{\|\eps\|_2}{\sigma})^2}{2 n}
    + (2 |T| + |Z\tilde T| ) 2 n^{-1/2}
    $.
    Using that $Z\sim N(0,1)$ and $\E[T^2]\le 1$, $\E[\tilde T^2]\le 1$,
    \begin{equation}
        \label{bound-rem-prim}
        \E[\Rem']
        \le 2\sqrt \mu + n^{-1/2}C
    \end{equation}
    for an absolute constant $C>0$.
    Combining \eqref{combined}-\eqref{combined2},
    $$
    \Bigl(1-\frac{\df}{n}\Bigr) + \Rem + \Rem'
    \ge \frac{\|y-X \hat b\|_2/\sqrt n}{(\sigma^2 + R)^{1/2}}
    -\sqrt \mu
    - \frac{\sqrt\mu R}{2\sigma^2}
    $$
    With the bounds in expectation
    on $\Rem$ in \eqref{bound-Rem} and on $\Rem'$
    in \eqref{bound-rem-prim},
    the lemma is proved.
\end{proof}

\bibliographystyle{plainnat}
\bibliography{../../bibliography/db}

\end{document}